\documentclass{amsart}

\usepackage{amssymb}
\usepackage{amsmath}
\usepackage{mathrsfs}
\usepackage{dsfont}
\usepackage{url}

\newtheorem{theorem}{Theorem}
\newtheorem*{theorem*}{Theorem}
\newtheorem{lemma}{Lemma}
\newtheorem*{lemma*}{Lemma}
\newtheorem{corollary}{Corollary}
\numberwithin{equation}{section}

\def \re {\operatorname{Re}}

\begin{document}

\title[The tail of the singular series]{The Tail of the Singular Series for the Prime Pair and Goldbach Problems}

\author[D. A. Goldston]{D. A. Goldston$^{*}$}

\address{Department of Mathematics,
 San Jos\'{e} State University,
 315 MacQuarrie Hall,
 One Washington Square,
 San Jos\'{e}, California 95192-0103,
 United States of America}
\date{\today}
\email{goldston@math.sjsu.edu}

\thanks{$^{*}$ Research supported by National Science Foundation Grants DMS-0804181 and DMS-1104434}
\author[Julian Ziegler Hunts]{Julian Ziegler Hunts}
\author[Timothy Ngotiaoco]{Timothy Ngotiaoco}



\begin{abstract} We obtain an asymptotic formula for a weighted sum of the square of the tail in the singular series for the Goldbach and prime-pair problems.
\end{abstract}

\maketitle

\thispagestyle{empty}

\section{Introduction and statement of results}

Hardy and Littlewood \cite{HardyLittlewood1922} conjectured in 1922 an asymptotic formula for the number of pairs of primes differing by $k$. The first major step forward on this conjecture only occurred in 2013 when Zhang \cite{Zhang2014} proved that there exist some $k$'s for which there are infinitely many such pairs of primes. Let  $\Lambda(n)$ be the von Mangoldt function defined by $\Lambda(n)=\log p$ if $n=p^m$, $p$ a prime, $m\ge 1$ an integer, and $\Lambda(n)=0$ otherwise. Hardy and Littlewood's conjecture is equivalent, for $k$ even, to
 
\begin{equation}\label{eq1.1} \psi_2(N,k):= \sum_{\substack{n,n^\prime\le N\\ n^\prime-n=k}}\Lambda(n)\Lambda(n^\prime) \sim \mathfrak{S}(k) (N-|k|) \quad \text{as} \ \ N\to \infty, \end{equation}
where 
\begin{equation}\label{eq1.2}
\mathfrak{S}(k) = \begin{cases}
{\displaystyle 2 C_{2} \prod_{\substack{ p \mid k \\ p > 2}} \!\left(\frac{p - 1}{p - 2}\right)} & \mbox{if $k$ is  even, $k\neq 0$,} \\
      0 & \mbox{if $k$ is  odd}\end{cases}
\end{equation}
and
\begin{equation}\label{eq1.3}
C_{2}
 = \prod_{p > 2}\! \left( 1 - \frac{1}{(p - 1)^{2}}\right)
 = 0.66016\ldots.
\end{equation}
For odd $k$ the sum in \eqref{eq1.1} has non-zero terms only when $n$ or $n^\prime$ is a power of 2, so $\psi_2(N,k) = O((\log N)^2).$  For the Goldbach problem Hardy and Littlewood conjectured an analogous formula for the number of ways an even number $k$ can be expressed as the sum of two primes, which also includes the arithmetic function $\mathfrak{S}(k)$.

The function $\mathfrak{S}(k)$ is called the singular series, a name given it by Hardy and Littlewood  because it first occurred as the series 
\begin{equation} \label{eq1.4} \mathfrak{S}(k) = \sum_{q=1}^\infty \frac{\mu(q)^2}{\phi(q)^2}c_q(-k),\end{equation}
where the Ramanujan sum $c_q(n)$ is defined by
\begin{equation} \label{eq1.5}   c_q(n) = \!\sum_{\substack{1\le a\le q\\ (a,q)=1}}\!\!\!\;e\Big(\frac{an}{q}\Big), \quad \ \  e(\alpha)=e^{2\pi i \alpha}. \end{equation} 
Some well-known properties of $c_q(n)$ (see, e.g., \cite{MontgomeryVaughan2007}) are that  $c_q(-n) = c_q(n)$,  $c_q(n)$
is a multiplicative function of $q$, and 
\begin{equation}\label{eq1.6} c_q(n) = \sum_{\substack{d|n\\d|q}} d\mu\!\left(\frac{q}{d}\right)\! = \frac{ \mu\big(\frac{q}{(n,q)}\big)\phi(q) }{\phi\big(\frac{q}{(n,q)}\big)} .\end{equation} 
Since the singular series is a sum of multiplicative functions in $q$, it is easy to verify that \eqref{eq1.4} is equivalent to the product in \eqref{eq1.2}. The series in \eqref{eq1.4} is a Ramanujan series; many arithmetic functions can be expanded into these series which have the property that the first term $q=1$ is the average or expected value of the arithmetic function. Thus we see that the $q=1$ term in \eqref{eq1.4} says that $\mathfrak{S}(k)$ has average value $1$. If we consider the first two terms we have
\[ \mathfrak{S}(k) =  1+ e\Big(\!\!\!\:-\frac{k}{2}\!\;\Big)\!\!\: +\sum_{q=3}^\infty \frac{\mu(q)^2}{\phi(q)^2}c_q(-k),\]
and therefore we obtain the refinement that on average $\mathfrak{S}(k)$ is $0$ if $k$ is odd and is $2$ if $k$ is even. 

In applications it is often useful to truncate the singular series; we write
\begin{equation}\label{eq1.7}\mathfrak{S}(k) = \mathfrak{S}_y(k)+\widetilde{\mathfrak{S}}_y(k), \end{equation}
where
\begin{equation} \label{eq1.8}  \mathfrak{S}_y(k) = \sum_{q\le y} \frac{\mu(q)^2}{\phi(q)^2}c_q(-k), \quad \widetilde{\mathfrak{S}}_y(k) = \sum_{q >y} \frac{\mu(q)^2}{\phi(q)^2}c_q(-k).\end{equation}
We refer to $\widetilde{\mathfrak{S}}_y(k)$ as the tail of the singular series. Montgomery and Vaughan \cite{MontgomeryVaughan1973}, by a simple argument using \eqref{eq1.5}, proved for $y\ge1$ the bound
\begin{equation} \label{eq1.9} 
\widetilde{\mathfrak{S}}_y(k) \ll  d(k) \frac{(\log\log 3y)^2}{y}.\end{equation}
Using a result of Ramanujan (for a proof see \cite{Wilson1922})
\[ \sum_{k\le N} d(k)^2 \sim \frac{1}{\pi^2} N(\log N)^3,\]
this bound immediately gives the mean square estimate 
\[ \sum_{k\le N} \widetilde{\mathfrak{S}}_y(k)^2 \ll \frac{N(\log N)^3 (\log\log 3y)^4}{y^2}. \]
In \cite{Goldston1992}\footnote{Beware that in \cite{Goldston1992} $\mathfrak{S}_y(k)$ and $\widetilde{\mathfrak{S}}_y(k)$ are defined differently than they are in this paper.} the first-named author improved this last bound by showing
\begin{equation} \label{eq1.10} \sum_{k\le N} \widetilde{\mathfrak{S}}_y(k)^2 \ll \frac{N\log N}{y^2}. \end{equation} 
Bounds of this type are useful in applications related to both the Goldbach and prime pair conjectures. For a recent application, see \cite{Buttkewitz}.
The proof of \eqref{eq1.10} is rather complicated and left open the question of whether the result can be improved further or is best possible. Our first result answers this question in the range $1\le y \le \sqrt{N}$. 
\begin{theorem} We have, for $1\le y\le \sqrt{N}$ and any fixed $\delta$, $0<\delta <1$, 
\begin{equation}\label{eq1.11} \sum_{k\le N}\left(N-k\right)^2 \widetilde{\mathfrak{S}}_y(k)^2 = \mathcal{T}(y) \frac{N^3}{3}\bigg(1+O_\delta\bigg(\!\!\:\Big(\frac{y^2}{N}{\Big)}\rule{0in}{.18in}^{\!\delta}\bigg)\bigg), \end{equation}
where
\begin{equation}\label{eq1.12} \mathcal{T}(y) := \sum_{q>y}\frac{\mu(q)^2}{\phi(q)^3}. \end{equation}
\end{theorem}
From \eqref{eq2.9} below we have
\begin{equation} \label{eq1.13} \mathcal{T}(y) = \frac{\mathcal{A}}{y^2}\left(1 +o(1)\right),\qquad \text{where}\ \  \mathcal{A} = \prod_p\left(1 +\frac{2-1/p}{(p-1)^2}\right)\! .\end{equation} 
A simple argument then gives the following result. Here $f \asymp g $ means $f\ll g$ and $g\ll f$.
\begin{corollary} We have, for some sufficiently small constant $c$,  
\begin{equation} \label{eq1.14}\sum_{k\le N} \widetilde{\mathfrak{S}}_y(k)^2 \asymp  \frac{N}{y^2}, \qquad \text{for} \quad 1\le y\le c\sqrt{N}.\end{equation}
and for $1\le y\le \sqrt{N}$ and any fixed $\delta$, $0<\delta <1$, 
\begin{equation}\label{eq1.15} \sum_{k\le N} \widetilde{\mathfrak{S}}_y(k)^2 = \mathcal{T}(y) N\bigg(1+O_\delta\bigg(\!\!\:\Big(\frac{y^2}{N}{\Big)}\rule{0in}{.18in}^{\!\delta/4}\bigg)\bigg). \end{equation}
\end{corollary}

Our main result is a refinement of Theorem 1.
\begin{theorem} We have, for $1\leq y\leq \sqrt{N}$,
\begin{equation} \label{eq1.16} \sum_{k\le N}\left(N-k\right)^2 \widetilde{\mathfrak{S}}_y(k)^2 = \mathcal{T}(y) \frac{N^3}{3} -\frac14 N^2 \!\left(\!\!\:\log\frac{N}{y^2}\right)^2\! + cN^2 \log \frac{N}{y^2} +O(N^2)+O\!\left(\frac{N^2}{\sqrt{y}}\log (2N)\!\!\:\right), \end{equation}
where 
\[ c = \frac34 -\frac12\log 2\pi + \frac12 \sum_p\frac {(p-2)\log p}{p(p-1)^2} .\]
\end{theorem}
The proof of Theorem 2 requires a less direct approach than Theorem 1. To proceed from the proof of Theorem 1 we want to take the parameter $\delta \ge 1$, but then the sums that result from the tail of the singular series diverge. Therefore we are forced to consider  $\widetilde{\mathfrak{S}}_y(k)^2=(\mathfrak{S}(k) - \mathfrak{S}_y(k))^2$, multiply this out, and evaluate each of the three terms separately. 

With a little additional work, by not dropping lower-order terms in \eqref{eq6.10.0}, \eqref{eq7.17.5}, \eqref{eq7.18} and \eqref{eq7.23} we can replace the $O(N^2)$ in Theorem 2 by $CN^2+O_\epsilon(N^2y^{-\frac12+\epsilon})$ for a complicated constant $C$.

The weight $(N-k)^2$ in our sum was chosen because it occurs naturally in the prime-pair problem. Obviously other weights or families of weights can be used. 

We have not been able to extend these results to the range $\sqrt{N}\le y \le N$ so in this range \eqref{eq1.10} remains the best result known.  For $y\ge N$, the method of \cite{Goldston1992} yields $\displaystyle\sum_{k\leq N}\widetilde{\mathfrak{S}}_y(k)^2\ll_A\frac{N\log N}{y^2\log(2y/N)^A}$.

{\it Notation.} We follow some common conventions. A sum will normally be over integers; any sum without a lower bound on the summation variable will start at 1. Empty sums will equal 0 and empty products will equal 1. The letter $p$ will always denote a prime. The letter $\epsilon$ will denote a small positive real number which may change from equation to equation.  

\section{Lemmas}
We gather here some of the results we need later.  
 \begin{lemma} Let $s(x)= x -\lfloor x\rfloor-\frac12$. Then for $x\ge 0$ we have 
\begin{equation} \label{eq2.1} \sum_{1\le k\le x}(x-k)= \frac12 \left(\!\!\;\Big(x-\frac12\Big)^{\!2}\!\!\; - s(x)^2\right)\end{equation}
and 
\begin{equation} \label{eq2.2}\sum_{1\le k\le x}(x-k)^2 = \frac13 \Big(x-\frac12\Big)^{\!3}\!\!\; - \int_{\frac12}^x s(u)^2\, du .\end{equation}
Since $|s(x)|\le \frac12$, we have
\begin{equation} \label{eq2.3} \sum_{1\le k\le x}(x-k)= \frac12 x^2-\frac12 x +O(1).\end{equation}
Since $s(x)$ is periodic with period $1$ and $\int_0^1 s(u)^2 \, du = \frac1{12}$, we have 
\begin{equation} \label{eq2.4} \sum_{1\le k\le x}(x-k)^2= \frac13 x^3 -\frac12x^2 +\frac16 x +O(1).\end{equation}
\end{lemma}

\noindent{\it Proof.} 
For the first identity, we use $\lfloor x\rfloor =x-\frac12 -s(x)$ to write
\[ \begin{split} S_1(x) := \sum_{1\le k\le x}\!\!\;(x-k) &= \sum_{1\le k\le \lfloor x\rfloor}\!\!\:(x-k)\\&= x\lfloor x\rfloor - \frac{\lfloor x\rfloor(\lfloor x\rfloor+1)}{2} \\ & 
= \frac12 \lfloor x\rfloor\left(\!\!\;x - \frac12 +s(x)\!\!\:\right)\\&
= \frac12 \left(\!\!\;x - \frac12 -s(x)\!\!\:\right)\!\left(\!\!\;x - \frac12 +s(x)\!\!\:\right)\\&= \frac12 \left( \!\!\:\Big(x-\frac12\Big)^{\!2}\!\!\; -s(x)^2\right).\end{split} \]
For the second identity, we use the first in
\[\sum_{1\le k\le x}\!\!\;(x-k)^2 = 2\int_{\frac12}^x\!\!\: S_1(u)\, du .\]

\begin{lemma} For fixed real numbers $a$ and $b$, let
\begin{equation}\label{eq2.5} G(x;a,b) := \sum_{r\le x}\frac{\mu(r)^2r^a}{\phi(r)^b}, \end{equation}
and 
\begin{equation} \label{eq2.6} g(s;a,b) := \prod_{p}\left( 1 - \frac{1 - p^{s-a+b}(1 - (1-\frac1p)^b)}{(p-1)^bp^{2(s-a)+b}}\right). \end{equation} 
Then we have
\begin{equation}\label{eq2.7} G(x;a,b) = \begin{cases}\frac{g(a-b+1;a,b)}{a-b+1}x^{a-b+1}+o_{a,b}(x^{a-b+1}) &\text{if }a-b>-1,\\ g(0;b-1,b)\log x +O_{a,b}(1)&\text{if }a= b-1,\\ \zeta(b-a) g(0;a,b)+ \frac{g(a-b+1;a,b)}{a-b+1}x^{a-b+1} + o_{a,b}(x^{a-b+1}) &\text{if }a-b <-1,\end{cases}\end{equation}
where $\zeta(s)$ is the Riemann zeta function \eqref{eq3.3}. 
\end{lemma}
\noindent This is Lemma 2 of \cite{Goldston1990}. In this paper we frequently apply this lemma to obtain only an upper bound for $G(x;a,b)$, but it is useful to know that the estimates obtained are essentially sharp. We note that when $a-b<-1$
\begin{equation}\label{eq2.9} \sum_{r> x}\frac{\mu(r)^2r^a}{\phi(r)^b} = \lim_{y\to\infty}\left(G(y;a,b)-G(x;a,b)\right)= \frac{g(a-b+1;a,b)}{b-a-1}x^{a-b+1} + o_{a,b}(x^{a-b+1}).\end{equation}

\begin{lemma}[Hildebrand] For $x\ge 1$, $d\geq1$,  we have 
\begin{equation} \label{eq2.10} \sum_{\substack{q\le x\\(q,d)=1}} \frac{\mu^2(q)}{\phi(q)} =  \frac{\phi(d)}{d}\Bigg(\!\!\!\;\log x  + \gamma +\sum_p\frac{\log p }{p(p-1)}+\sum_{p\vert d}\frac{\log p}{p}\;\!\Bigg)+ O\Big( x^{-\frac12}\prod_{p\vert d}(1+p^{-\frac12})\Big). \end{equation}
\end{lemma}
\noindent This is Hilfssatz 2 of \cite{Hildebrand1984}.
\begin{lemma} For $x\ge 1$,
\begin{equation} \label{eq2.11} \sum_{k\le x} (x-k)\mathfrak{S}(k)= \frac12 x^2 -\frac12 x\log x + 
\frac12(1-\gamma -\log 2\pi)x + O_\epsilon(x^{\frac12+\epsilon}).\end{equation}
\end{lemma} 
\noindent This was first stated in \cite{Goldston1990}, and also appeared in \cite{FriedlanderGoldston1995}, but the first published proof is in \cite{MontgomerySound1999}.

Our next lemma is a generalization and strengthening of Lemma 4 due to Vaughan. We let 
\begin{equation} \label{eq2.12}  \mathfrak{G}_d(k) = 2C(d)\!\prod_{\substack{ p \mid k \\ (p,2d)=1}}\!\!\!\!\: \left(\frac{p - 1}{p - 2}\right), \end{equation}
where 
\begin{equation} \label{eq2.13} C(d) =\! \prod_{(p,2d)=1}\!\!\!\!\: \left( 1 - \frac{1}{(p - 1)^{2}}\right).\end{equation}
Note that unlike for $\mathfrak{S}(k)$ we do not require that $\mathfrak{G}_d(k)$ be zero if $k$ is odd; instead $\mathfrak{G}_d(k)=\mathfrak{G}_d(2k)$.
\begin{lemma}[Vaughan] For $x\ge 1$, we have 
\begin{equation} \label{eq2.14} \sum_{k\le x} (x-k)\mathfrak{G}_d(k) = x^2 - \frac12 \frac{(d,2)\phi(d)}{d}x \Bigg(\!\!\!\; \log x  +\gamma -1+\log 2\pi+ \sum_{p|2d}\frac{\log p}{p-1}\;\!\Bigg)\!\!\: + E(x,d)\end{equation}
where 
\begin{equation} \label{eq2.15} E(x,d) \ll x^{\frac12}\exp\!\!\;\left(-c\frac{(\log 2x)^{\frac35}}{(\log\log 3x)^{\frac15}}\right)\prod_{p|d}(1- p^{-\frac14})^{-1}\end{equation}
for some positive constant $c$.  If we assume the Riemann Hypothesis then $x^{\frac12}$ in \eqref{eq2.15} can be replaced by $x^{\frac5{12}+\epsilon}$.
\end{lemma}
\noindent This is Theorem 3 of \cite{Vaughan2001}. (The Riemann Hypothesis estimate is on page 552 of that paper.) We can recover Lemma 4 from Lemma 5 with a stronger error term by using
\[\sum_{k\le x} (x-k)\mathfrak{S}(k)= 2\sum_{k\le \frac{x}{2}} \!\Big(\frac{x}{2}-k\Big)\mathfrak{G}_1(k).\]

\section{Proof of Theorem 1.} 
We have
\begin{equation}\label{3.1} S := \sum_{k\le N}(N-k)^2\widetilde{\mathfrak{S}}_y(k)^2 =\sum_{q> y}\sum_{q^\prime> y}\frac{\mu(q)^2}{\phi(q)^2}\frac{\mu(q^\prime)^2}{\phi(q^\prime)^2}\!\!\;\underbrace{\sum_{1\leq k\leq N}\!\!\;(N-k)^2c_q(-k)c_{q^\prime}(-k)}_{S^\prime},\end{equation}
and by the formula $c_q(-k)=\displaystyle\sum_{\substack{d\mid q\\d\mid k}}d\!\;\mu\!\left(\frac{q}{d}\right)$, we have 
\begin{equation}\label{eq3.2}S^\prime=\sum_{d\mid q}\sum_{d^\prime\mid q^\prime}d\!\;\mu\!\left(\frac{q}{d}\right)\!\!\;d^\prime\mu\Big(\frac{q^\prime}{d^\prime}\!\!\;\Big)\!\sum_{\substack{1\leq k\leq N\\ [d,d^\prime]\mid k}}(N-k)^2.\end{equation} 
We now need to evaluate the inner sum over $k$. In proving Theorem 2 we do this with the elementary  Lemma 1, but here we need to use the formula in Theorem B of Ingham \cite{Ingham1932}: if  $m$ is a positive integer,  $c>0$, and $x>0$, then 
\[  \frac{m!}{2\pi i}\int_{c-i\infty}^{c+i\infty} \frac{x^{s+m}}{s(s+1)(s+2)\cdots (s+m)}\, ds  = \begin{cases}0 &\text{if } 0<x\le 1, \\ (x-1)^m &\text{if }x\ge 1.\end{cases}  \]
The Riemann zeta function is, for $s=\sigma +it$, $\sigma>1$,
\begin{equation} \label{eq3.3} \zeta(s) = \sum_{n=1}^\infty \frac{1}{n^s} = \prod_p\left(1-\frac1{p^s}\right)^{\!-1}\!\!.\end{equation}
The series and product converge absolutely and converge uniformly for $\sigma\ge1+\epsilon$. Hence for $x\ge 1$ and $c>1$ we have
\begin{equation} \label{eq3.4} \sum_{1\le n\le x} (x-n)^k = \frac{k!}{2\pi i}\int_{c-i\infty}^{c+i\infty}\!\!\: \frac{\zeta(s)x^{s+k}}{s(s+1)(s+2)\cdots (s+k)}\, ds. \end{equation}
Now
\[ \sum_{\substack{1\leq k\leq N\\ [d,d^\prime]\mid k}}\!\!\:(N-k)^2=[d,d^\prime]^2\!\!\sum_{\substack{1\leq k\leq \frac{N}{[d,d^\prime]}}}\!\!\!\:\left(\frac{N}{[d,d^\prime]}-k\right)^2\!\!\:,\]
and therefore
\[ \sum_{\substack{1\leq k\leq N\\ [d,d^\prime]\mid k}}\!\!\:(N-k)^2= \frac{2!}{2\pi i}\int_{c-i\infty}^{c+i\infty}\!\!\: \frac{\zeta(s)N^{s+2}}{s(s+1)(s+2)[d,d^\prime]^s}\, ds, \]
making use of the assumption that $y\le \sqrt{N}$ to ensure that
\[\frac{N}{[d,d^\prime]}\ge \frac{N}{dd^\prime}\ge \frac{N}{qq^\prime}\ge \frac{N}{y^2}\ge 1.\]
Hence 
\begin{equation} \label{eq3.5} S = \frac{2!}{2\pi i}\int_{c-i\infty}^{c+i\infty}\!\!\: \zeta(s) B_s(y)\frac{N^{s+2}}{s(s+1)(s+2)}\, ds,  \end{equation}
where 
\begin{equation} \label{eq3.6} B_s(y) =\sum_{q> y}\sum_{q^\prime> y}\frac{\mu(q)^2}{\phi(q)^2}\frac{\mu(q^\prime)^2}{\phi(q^\prime)^2}\sum_{d\mid q}\sum_{d^\prime\mid q^\prime}\frac{dd^\prime \mu\!\left(\frac{q}{d}\right)\!\!\:\mu\big(\frac{q^\prime}{d^\prime}\!\!\;\big)}{[d,d^\prime]^s}.\end{equation}
Following the method Selberg introduced for the Selberg sieve \cite{Selberg1947}, we now diagonalize $B_s(y)$. Define $\phi_s(n)$ by the equation $n^s = \sum_{d\mid n} \phi_s(d)$, so that $\phi_s(n)$ is multiplicative and $\phi_s(p) = p^s-1$. Letting $n=(d,d^\prime)$ we have 
\[ (d,d^\prime)^s= \sum_{\substack{r\mid d\\ r\mid d^\prime}} \phi_s(r),\]
and thus
\[ \frac{dd^\prime}{{[d,d^\prime]}^s} = (dd^\prime)^{1-s}\sum_{\substack{r\mid d\\ r\mid d^\prime}} \phi_s(r). \]
Hence
\begin{equation}\label{eq3.7} B_s(y) =\sum_{r=1}^\infty \phi_s(r)\!\!\:\left(\rule{0in}{.35in}\right.\! \sum_{\substack{q> y\\ r\mid q}}\frac{\mu(q)^2}{\phi(q)^2}\sum_{\substack{d\mid q\\ r\mid d}}d^{1-s}\mu\!\left(\frac{q}{d}\right)\!\!\!\!\;\left.\rule{0in}{.35in}\right)^{\!\!\!\;2}\!.\end{equation}

The simplest bound on $\zeta(s)$ in the critical strip is that if $0<\alpha < 1$, $|t|\ge 1$, then  
\begin{equation} \label{eq3.8} |\zeta(s)| < C(\alpha) |t|^{1-\alpha}\quad\text{for }\sigma\ge\alpha\end{equation} 
for some constant $C(\alpha)$, see Theorem 9 of Ingham\cite{Ingham1932}. We also need the bound, for $0<\alpha < 1$,  
\begin{equation}\label{eq3.9} B_{s}(y)\ll  \frac1{y^{2\alpha}}\ \ \text{for} \  \ \sigma \ge \alpha ,\end{equation} 
which we will prove later in this section. 
In our formula for $S$ we move the contour to the left past the simple pole with residue 1 at $s=1$ of $\zeta(s)$  to the line $s=\alpha +it$ with $0<\alpha <1$. Since by \eqref{eq3.8} and \eqref{eq3.9} the integrand is $O_{\alpha}(N^{2+\alpha}y^{-2\alpha}/|t|^{2+\alpha})$ for $|t|\ge1$, the integrals converge absolutely and we obtain
\begin{equation}\begin{split} \label{eq3.10} S & = B_1(y)\frac{N^3}{3}  + \frac{2!}{2\pi i}\int_{
\alpha -i\infty}^{\alpha +i\infty}\!\!\: \zeta(s) B_s(y)\frac{N^{s+2}}{s(s+1)(s+2)}\, ds \\ &= B_1(y)\frac{N^3}{3}+O_\alpha\!\left(\!\!\;\frac{N^{2+\alpha}}{y^{2\alpha}}\!\!\;\right)\!. \end{split}  \end{equation}
We have $\displaystyle\sum_{\substack{d\mid q\\r\mid d}}\mu\!\left(\frac{q}{d}\right)=\sum_{s\mid \frac{q}{r}}\mu\Big(\frac{q/r}{s}\Big)=\begin{cases}1&\text{if }q=r,\\0&\text{if }q\neq r,\end{cases}$ and thus
\begin{equation}\label{eq3.12}B_1(y)=\sum_{r=1}^\infty \phi(r)\!\!\:\left(\rule{0in}{.35in}\right.\!\sum_{\substack{q> y\\r\mid q}}\frac{\mu(q)^2}{\phi(q)^2}\sum_{\substack{d\mid q\\r\mid d}}\mu\!\left(\frac{q}{d}\right)\!\!\!\!\;\left.\rule{0in}{.35in}\right)^{\!\!\!\;2}=\sum_{r=1}^\infty \phi(r)\!\!\:\left(\rule{0in}{.35in}\right.\!\sum_{\substack{q> y\\\vphantom{\mid}q=r}}\frac{\mu(q)^2}{\phi(q)^2}\!\!\!\!\;\left.\rule{0in}{.35in}\right)^{\!\!\!\;2}=\sum_{r> y}\frac{\mu(r)^2}{\phi(r)^3}= \mathcal{T}(y).\end{equation}
We conclude that 
\[ S = \mathcal{T}(y) \frac{N^3}{3} + O_{\alpha}\!\left(\frac{N^3}{y^2} \left(\frac{y^2}{N}\right)^{1-\alpha}\right)\!,\]
which proves Theorem 1 on taking $\alpha = 1-\delta$.

It remains to prove \eqref{eq3.9}.  For the sums over $q$ and $d$ inside the square in \eqref{eq3.7},  writing $q=du$, $d=rv$, we have $q=ruv$ and 
\[\sum_{\substack{q> y\\ r\mid q}}\frac{\mu(q)^2}{\phi(q)^2}\sum_{\substack{d\mid q\\ r\mid d}}d^{1-s}\mu\!\left(\frac{q}{d}\right)= \sum_{ruv> y} \frac{\mu(ruv)^2}{\phi(ruv)^2}(rv)^{1-s}\mu(u).\]
Hence
\[ B_s(y) = \sum_{r=1}^\infty \frac{\mu(r)^2\phi_s(r)r^{2-2s}}{\phi(r)^4}\!\!\;\left(\rule{0in}{.35in}\right.\!\sum_{\substack{u=1 \\ (u,r)=1}}^\infty \frac{\mu(u)}{\phi(u)^2} \left(\rule{0in}{.35in}\right.\!\sum_{\substack{v> y/ur \\ (v,ur)=1}} \frac{\mu(v)^2v^{1-s}}{\phi(v)^2}\!\!\!\!\;\left.\rule{0in}{.35in}\right)\!\!\!\!\;\left.\rule{0in}{.35in}\right)^{\!\!\!\;2}\!.\]
We note that for squarefree $r$ 
\[ |\phi_{a+it}(r)|\le \prod_{p| r} \left(p^a +1\right)=r^a \prod_{p|r}\left( 1+ \frac{1}{p^a} \right)=r^a \sigma_{-a}(r),\]
where $\sigma_z(r)=\sum_{d|r}d^z$. 
We conclude
\[|B_{a+it}(y)|\le \sum_{r=1}^\infty \frac{\mu(r)^2r^{2-a}}{\phi(r)^4}\sigma_{-a}(r)\!\!\:\left(\sum_{u=1}^\infty \frac{\mu(u)^2}{\phi(u)^2} \!\!\;\left(\sum_{v> y/ur }\!\frac{\mu(v)^2v^{1-a}}{\phi(v)^2}\right)\!\right)^{\!2}\!.\]
Clearly the right-hand side is a decreasing function of $a$, and therefore to prove \eqref{eq3.9} we only need to prove that the right-hand side above satisfies the bound in \eqref{eq3.9} for $a=\alpha$. 
Since by Lemma 2 
\[ \sum_{v>y/ur} \!\frac{\mu(v)^2v^{1-\alpha}}{\phi(v)^2} \ll_\alpha\!\!\:\left(\frac{ur}{y}\right)^{\!\alpha}\!,\]
we have 
\[ B_{\alpha+it}(y)\ll_\alpha \sum_{r=1}^\infty \frac{\mu(r)^2 r^{2-\alpha}}{\phi(r)^4}\sigma_{-\alpha}(r)\!\!\:\left(\sum_{u=1}^\infty \frac{\mu(u)^2}{\phi(u)^2}\left(\frac{ur}{y}\right)^{\alpha} \right)^{\!2}\!.\]
Applying Lemma 2 again, the sum over $u$ is
 \[\frac{r^\alpha}{y^\alpha}\sum_{u=1}^\infty\frac{\mu(u)^2 u^\alpha}{\phi(u)^2}\ll_\alpha\frac{r^\alpha}{y^\alpha},\] 
so
\[B_s(y) \ll_\alpha\sum_{r=1}^\infty \frac{\mu(r)^2r^{2-\alpha} \sigma_{-\alpha}(r)}{\phi(r)^4}\frac{r^{2\alpha}}{y^{2\alpha}}=\frac{1}{y^{2\alpha}}\sum_{r=1}^\infty \frac{\mu(r)^2r^{2+\alpha} \sigma_{-\alpha}(r)}{\phi(r)^4}.\]
Since $\phi(dm)=\phi(d)\phi(m)$ when $\mu(dm)\neq0$,
\[\sum_{r=1}^\infty \frac{\mu(r)^2r^{2+\alpha}}{\phi(r)^4}\sum_{d|r} \frac1{d^\alpha}= \sum_{d,m} \frac{\mu(dm)^2d^2m^{2+\alpha}}{\phi(dm)^4 }\le\left(\sum_{m=1}^\infty\frac{\mu(m)^2m^{2+\alpha}}{\phi(m)^4}\right)\left( \sum_{d=1}^\infty \frac{\mu(d)^2d^2}{\phi(d)^4}\right)\ll_\alpha 1\]
and $B_s(y)\ll_\alpha\frac{1}{y^{2\alpha}}$, which proves \eqref{eq3.9}.

To prove Corollary 1, let 
\[  T_0(N) = \sum_{k\le N} \widetilde{\mathfrak{S}}_y(k)^2  , \qquad   T_m(N) = \sum_{k\le N}(N-k)^m\widetilde{\mathfrak{S}}_y(k)^2 \quad \text{for} \ \  m\ge 1 .\]
 Then by Theorem 1 and \eqref{eq1.13} for $1\le y\le cN^{1/2}$ with $c$ sufficiently small  
\[  T_2(N)  \asymp  \frac{N^3}{y^2} \]
and \eqref{eq1.14} follows from
\[  \frac{1}{N^2}T_2(N) \le T_0(N) \le \frac{1}{N^2}T_2(2N).\]
To prove \eqref{eq1.15} we note for $m\ge 0$ that
\[  \quad T_{m+1}(N)= (m+1)\int_1^N T_m(u)\, du  .\]
Since $T_m(N)$ is a nondecreasing function of $N$,
we have, for $1\le h \le N$,  
\[ T_1(N) \le \frac{1}{h} \int_{N}^{N+h} T_1(u)\, du = \frac{T_2(N+h) -T_2(N)}{2h} \]
and similarly 
\[ T_1(N)\ge \frac{T_2(N) -T_2(N-h)}{2h}.\]
Now by \eqref{eq1.11} and \eqref{eq1.13}
\[ T_2(N) = \mathcal{T}(y) \frac{N^3}{3} + O_\delta\!\left( \frac{N^3}{y^2}\left(\frac{y^2}{N}\right)^\delta\right)\!,\]
and hence
\[ \begin{split}  \frac{T_2(N \pm h) - T_2(N)}{\pm 2h} &=\frac12 \mathcal{T}(y)\left( N^2 \pm N h + \frac{h^2}{3}\right) + O_\delta\!\left( \frac{N^3}{hy^2}\left(\frac{y^2}{N}\right)^\delta\right)\\& = \mathcal{T}(y) \frac{N^2}{2} + O\!\left(\frac{Nh}{y^2}\right)\! + O_\delta\!\left( \frac{N^3}{hy^2}\left(\frac{y^2}{N}\right)^\delta\right)\!.
\end{split} \] 
Balancing the two error terms by choosing $h = N \left(\frac{y^2}{N}\right)^{
\frac{\delta}{2}}$, we conclude
\[  \frac{T_2(N \pm h) - T_2(N)}{\pm 2h} =   \mathcal{T}(y) \frac{N^2}{2} + O_\delta\!\left( \frac{N^2}{y^2}\left(\frac{y^2}{N}\right)^{\frac{\delta}{2}}\right),\]
and hence 
\[ T_1(N) =   \mathcal{T}(y)\frac{N^2}{2} + O_\delta\!\left( \frac{N^2}{y^2}\left(\frac{y^2}{N}\right)^{\frac{\delta}{2}}\right).\]
By the same argument $T_0(N)$ is bounded between the expressions 
\[ \frac{T_1(N \pm h) - T_1(N)}{\pm h} =  \mathcal{T}(y)N + O\!\left(\frac{h}{y^2}\right)\!+  O_\delta\!\left( \frac{N^2}{hy^2}\left(\frac{y^2}{N}\right)^{\frac{\delta}{2}}\right)\] 
and the choice $h = N \left(\frac{y^2}{N}\right)^{
\frac{\delta}{4}}$ gives 
\[ T_0(N) = \mathcal{T}(y)N + O_\delta\!\left(\frac{N}{y^2} \left(\frac{y^2}{N}\right)^{
\frac{\delta}{4}}\right),\]
which proves \eqref{eq1.15}.

\section{ The average of the singular series tail}

In this section for completeness we give a proof of the average size of the tail of the singular series. This proof illustrates the method we use to prove Theorem 2 without all the complications.
\begin{theorem} We have, for $1\le y\le N$,  
\begin{equation} \label{eq4.1} \sum_{k\le N} (N-k)\widetilde{\mathfrak{S}}_y(k) = -\frac12 N\log \frac{N}{y} + \frac12\left(1 -\log 2\pi +\sum_p \frac{\log p}{p(p-1)}\right) N +O(Ny^{-\frac12}) +O(y).\end{equation}
\end{theorem}

The reason the average does not have a main term of size $\frac{N^2}{y}$ as one might expect is that the term 
$1$ from $q=1$ in \eqref{eq1.4} cancels out this term independent of the truncation level $y$.

{\it Proof.} We have 
\begin{equation} \label{eq4.2} \sum_{k\le N} (N-k)\widetilde{\mathfrak{S}}_y(k)= 
\sum_{k\le N} (N-k)\mathfrak{S}(k) - \sum_{k\le N} (N-k)\mathfrak{S}_y(k).\end{equation}
The first sum is evaluated in Lemma 4. 
For the second sum, we use \eqref{eq1.8} and \eqref{eq1.6} to obtain 
\[ \begin{split}\sum_{k\le N} (N-k)\mathfrak{S}_y(k) &= \sum_{q\le y}\frac{\mu(q)^2}{\phi(q)^2} \sum_{k\le N} (N-k)c_q(-k)\\ &= \sum_{q\le y}\frac{\mu(q)^2}{\phi(q)^2}\sum_{d|q} d \mu\!\left(\frac{q}{d}\right)\bigg( \sum_{\substack{1\le k\le N \\ d|k}} (N-k) \bigg), \end{split}\]
and we see on letting $k=dm$ that by Lemma 1  
\[ \begin{split}\sum_{\substack{1\le k\le N \\ d|k}} (N-k) &= d\sum_{ 1\le m\le \frac{N}{d} } \left( \frac{N}{d} -m\right) \\&
=\frac12 \frac{N^2}{d} -\frac12N +O(d), \end{split}\]
and hence
\[ \sum_{k\le N} (N-k)\mathfrak{S}_y(k) = \frac12 N^2\sum_{q\le y}\frac{\mu(q)^2}{\phi(q)^2}\sum_{d|q}\mu\!\left(\frac{q}{d}\right) -\frac12 N \sum_{q\le y}\frac{\mu(q)^2}{\phi(q)^2}\sum_{d|q}d\mu\!\left(\frac{q}{d}\right)+O\left(\sum_{q\le y}\frac{\mu(q)^2}{\phi(q)^2}\sum_{d|q}d^2\right) .\]
By Lemma 2,
\[\begin{split} \sum_{q\le y}\frac{\mu(q)^2}{\phi(q)^2}\sum_{d|q}d^2 &= \sum_{dm\le y} \frac{\mu(dm)^2d^2}{\phi(dm)^2} \le\left(\sum_{m\le y}\frac{\mu(m)^2}{\phi(m)^2}\right)\left(\sum_{d\le y} \frac{\mu(d)^2 d^2}{\phi(d)^2}\right)\ll y. \end{split}  \] 
Hence we see
\begin{equation} \label{eq4.4} \sum_{k\le N} (N-k)\mathfrak{S}_y(k) = \frac12 N^2 -\frac12 N \sum_{q\le y}\frac{\mu(q)^2}{\phi(q)}+O(y). \end{equation}
The theorem now follows from \eqref{eq4.2}, \eqref{eq4.4}, Lemma 3 with $d=1$, Lemma 4, and the fact that $N^{\frac12+\epsilon}\le\max(Ny^{-\frac12},y)$ for $\epsilon\leq\frac16$. 

\section{Starting the proof of Theorem 2}  To prove Theorem 2 we need to asymptotically evaluate
\begin{equation}\label{eq5.1}\begin{split}  \sum_{k\le N}\left(N-k\right)^2 \widetilde{\mathfrak{S}}_y(k)^2 &=
\sum_{k\le N}\left(N-k\right)^2 \big(\mathfrak{S}(k)-\mathfrak{S}_y(k)\big)^2\\ &= \sum_{k\le N}\left(N-k\right)^2\mathfrak{S}(k)^2 -2 \sum_{k\le N}\left(N-k\right)^2\mathfrak{S}(k)\mathfrak{S}_y(k) + \sum_{k\le N}\left(N-k\right)^2\mathfrak{S}_y(k)^2\\ &
= : S_1 -2S_2 + S_3 . \end{split}\end{equation}
We evaluate each of these terms in the following sections.

\section{The sum $S_1$}In this section we evaluate 
\begin{equation} \label{eq6.1} S_1 = \sum_{k\le N}\left(N-k\right)^2\mathfrak{S}(k)^2 .
\end{equation}
The proof is along the same lines as the proof in \cite{MontgomerySound1999} of Lemma 4. 
\begin{theorem} We have 
 \begin{equation}\label{eq6.2}\begin{split} \sum_{k\le N}\left(N-k\right)^2\mathfrak{S}(k)^2 &= \prod_p\left( 1 + \frac{1}{(p-1)^3}\right) \frac{N^3}{3} -\frac14 N^2(\log N)^2 \\ & \ \ + \left(\frac34 -\gamma -\frac12\log 2\pi -\frac12\sum_p\frac{\log p}{(p-1)^2}\right)N^2\log N + O(N^2).\end{split} \end{equation}
\end{theorem}

\noindent {\it Proof.} Let $g(k) = \displaystyle\prod_{\substack{p|k\\p>2}} \left(\frac{p-1}{p-2}\right)^{\!2}$, so that $\mathfrak{S}(k)^2 = \begin{cases}4C_2^2 g(k) & \mbox{if } 2 \mid k,\\ 0 & \mbox{if } 2 \nmid k\end{cases}$ and
\begin{equation}\label{eq6.3}
S_1 = 4C_2^2\sum_{\substack{1 \le k \le N \\ 2\mid k}}(N-k)^2  g(k) = 16C_2^2\sum_{1 \le k \le \frac{N}{2}} \left(\frac{N}{2}-k\right)^{\!2}\!  g(k)
= 16C_2^2S_{11}\:\!\!\!\left(\frac{N}{2}\right)
\end{equation}
where $S_{11}(N) = \displaystyle{\sum_{1 \le k \le N} (N-k)^2 g(k)}.$  Let 
\begin{equation}\label{eq6.4}\begin{split}
G(s) = \sum_{n=1}^\infty \frac{g(n)}{n^s} &= \prod_p\left(1+\sum_{m=1}^\infty \frac{g(p^m)}{p^{ms}}\right)\\& = \left(1-\frac{1}{2^{s}}\right)^{-1}\prod_{p>2}\left(1+\left(\frac{p-1}{p-2}\right)^2\frac{1}{p^s-1}\right),
\end{split}\end{equation}
for $\re s>1$. To analytically continue $G(s)$ to the left, we see the dominant factor is  
\[ \zeta(s)= \prod_p\left(1- \frac{1}{p^s}\right)^{-1},\]
 and therefore we have
\begin{equation}\label{eq6.5}\begin{split}
G(s) &= \zeta (s) \prod_{p>2} \left(1+\left(\frac{p-1}{p-2}\right)^2 \frac{1}{p^s-1}\right)\left(1-\frac1{p^{s}}\right)\\
&= \zeta (s) \prod_{p>2} \left(1+\frac{2p-3}{(p-2)^2p^s}\right) =:\zeta (s) H(s)
\end{split}\end{equation}
with $H(s)$ analytic for $\re s>0$.  Next we write $\displaystyle H(s)=\zeta(s+1)^2\left(1-\frac{1}{2^{s+1}}\right)^2 J(s)$,
\begin{equation}\label{eq6.5.5} J(s)=\prod_{p>2} \left(1+\frac{2p-3}{(p-2)^2 p^s}\right)\left(1-\frac{1}{p^{s+1}}\right)^2.\end{equation}
We then have (by Mathematica) that
\begin{equation}\label{eq6.6}
J(s) =\prod_{p>2} \left(1+\frac{1}{(p-2)^2 p^s} \left(5 - \frac{8}{p} + \frac{1}{p^s} \left(-3 + \frac{2}{p} + \frac{4}{p^2}\right) + \frac{1}{p^{2s}} \left(\frac{2}{p}-\frac{3}{p^2}\right)\right)\right),
\end{equation}
and for $-1<\re s<0$ this is $\displaystyle{\prod_{p>2} \left(1+\frac{1}{(p-2)^2 p^s} \left(-\frac{3}{p^s}+O(1)+O\Big(\frac1{p^{2s+1}}\Big)\right)\right)}$, which is analytic for $\re s>-\frac{1}{2}$.  

Now, in the same way we obtained \eqref{eq3.4}, we have for $a>1$
\begin{equation}\label{eq6.7}
S_{11}(N)= \sum_{k=1}^N (N-k)^2 g(k) = \frac{2!}{2 \pi i} \int_{a - i\infty}^{a+ i\infty} G(s) \frac{N^{s+2}}{s(s+1)(s+2)}\,ds.
\end{equation}
We move the contour to $\re s=b$, $-\frac{1}{2} < b < 0$.  To ensure convergence and justify moving the contour we need to use a standard bound for $\zeta(s)$ which improves on \eqref{eq3.8}. By \cite{Titchmarsh}, Chapter 5, for $|t|\ge 1$, 
\begin{equation} \label{eq6.75}\zeta(\sigma+it) \ll (|t|+3)^{\lambda(\sigma)+\epsilon},\end{equation}
where
\begin{equation} \label{eq6.76}\lambda(\sigma)= \begin{cases}0 &\text{if }\sigma>1,\\ \frac12-\frac12\sigma &\text{if }0<\sigma \le 1, \\ \frac12-\sigma  &\text{if }\sigma \le 0.\end{cases}\end{equation}
This, along with the fact that $J=O_b(1)$ for $\re s\ge b$, shows that the integrand is $O_b(N^{\sigma+2}/|t|^{2b+\frac52})$ for $|t|\geq1$.

We encounter a simple pole at $s=1$ and a triple pole at $s=0$.  Since $H(s)$ is analytic at $s=1$ and $\zeta(s) = \dfrac{1}{s-1} + O(1)$, the pole at $1$ contributes $\frac{1}{3}H(1)N^3$ to $S_{11}(N)$.  Expanding around $s=0$ we have
\begin{equation}\label{eq6.8}\begin{split}
G(s)\frac{N^{s+2}}{s(s+1)(s+2)}&=\zeta(s)\zeta(s+1)^2\left(1-\frac{1}{2^{s+1}}\right)^2J(s)\frac{N^{s+2}}{s(s+1)(s+2)}\\&=\frac{N^2}{4}\cdot\frac{1}{s^3}K(s)N^s\\&=\frac{N^2}{4}\cdot\frac{1}{s^3}\left(1+(\log N)s+\frac{(\log N)^2}{2}s^2+O(s^3)\right)\\&\quad\;\;\;\;\cdot\left(K(0)+K^\prime(0)s+\frac{K^{\prime\prime}(0)}{2}s^2+O(s^3)\right)\!,
\end{split}\end{equation}
where
\[K(s)=\zeta(s)(s\zeta(s+1))^2(2-2^{-s})^2\frac{1}{(1+s)(2+s)}J(s).\]
The pole at 0 therefore contributes $\frac{N^2}{2}\left(\frac12K(0)(\log N)^2+K^\prime(0)\log N+\frac12K^{\prime\prime}(0)\right)$.  From the expansion $s\zeta(s+1)=1+\gamma s+O(s^2)$, we find that $K(0)=\frac12\zeta(0)J(0)=-\frac14J(0)$ and, using that if $f_1$ and $f_2$ are differentiable then $\frac{(f_1f_2)^\prime}{f_1f_2}=\frac{f_1^\prime}{f_1}+\frac{f_2^\prime}{f_2}$, that $\frac{K^\prime(0)}{K(0)}=\frac{\zeta^\prime(0)}{\zeta(0)}+2\gamma+2\log 2-1-\frac12+\frac{J^\prime(0)}{J(0)}$.  
We have
\begin{equation}\label{eq6.9}\begin{split}
J(0) &= \prod_{p>2} \left(1 + \frac{2p-3}{(p-2)^2}\right)\left(1-\frac{1}{p}\right)^2 = \prod_{p>2} \frac{((p-2)^2 + 2p - 3)(p-1)^2}{p^2(p-2)^2}\\&=\prod_{p>2} \frac{(p-1)^4}{p^2(p-2)^2} = \prod_{p>2} \left(\frac{1}{1-\frac{1}{(p-1)^2}}\right)^2 = \frac{1}{C_2^2},
\end{split}\end{equation}
$\frac{\zeta^{\prime}(0)}{\zeta(0)}=\log2\pi$,
\begin{equation}\label{eq6.9.5}\frac{J^\prime(0)}{J(0)}=\sum_{p>2}\left.\left(2\frac{(\log p)p^{-s}}{p-p^{-s}}-\frac{(2p-3)(\log p)p^{-s}}{(p-2)^2+(2p-3)p^{-s}}\right)\right\vert_{s=0}=\sum_{p>2}\frac{\log p}{(p-1)^2},\end{equation}
and 
\begin{equation}\label{eq6.9.7}\begin{split}H(1)&=\prod_{p>2}\left(1+\frac{2p-3}{(p-2)^2p}\right)=\prod_{p>2}\frac{p^3-4p^2+6p-3}{(p-2)^2p}\\&=\prod_{p>2}\frac{(p-1)^4+p-1}{p^2(p-2)^2}=\prod_{p>2}\frac{(p-1)^4}{p^2(p-2)^2}\left(1+\frac{1}{(p-1)^3}\right)\\&=\frac{1}{2C_2^2}\prod_p\left(1+\frac{1}{(p-1)^3}\right).\end{split}\end{equation}

Combining these, we obtain
\begin{equation}\label{eq6.10.0}\begin{split}S_1(N)&=16C_2^2S_{11}\:\!\!\!\left(\frac{N}{2}\right)\\&=16C_2^2\left(\frac{1}{24}H(1)N^3+\frac{N^2}{16}\left(K(0)(\log(N/2))^2+2K^\prime(0)\log(N/2)+K^{\prime\prime}(0)\right)\right)\\&\quad+\frac{2!}{2\pi i}\int_{b-i\infty}^{b+i\infty}G(s)\frac{N^{s+2}}{s(s+1)(s+2)}\,ds\\&=\frac{1}{3}\prod_p\left(1+\frac{1}{(p-1)^3}\right)N^3-\frac{1}{4}C_2^2J(0)N^2((\log N)^2-2(\log 2)\log N+(\log 2)^2)\\&\quad-\frac{1}{2}C_2^2J(0)\left(\log 2\pi+2\gamma+2\log 2-1-\frac{1}{2}+\sum_{p>2}\frac{\log p}{(p-1)^2}\right)N^2(\log N-\log 2)+O(N^2)\\&=\prod_p\left(1+\frac{1}{(p-1)^3}\right)\frac{N^3}{3}-\frac{1}{4}N^2(\log N)^2\\&\quad+\left(\frac{3}{4}-\gamma-\frac{1}{2}\log 2\pi-\frac{1}{2}\sum_p\frac{\log p}{(p-1)^2}\right)N^2\log N+O(N^2),\end{split}\end{equation}
as desired.

\section{The sum $S_2$}
In this section we evaluate
\begin{equation}\label{7.1} S_2=\sum_{k\leq N}(N-k)^2\mathfrak{S}(k)\mathfrak{S}_y(k).\end{equation}
\begin{theorem} We have
\begin{equation}\label{eq7.26}\begin{split}
S_2&=\left(\sum_{q\leq y}\frac{\mu(q)^2}{\phi(q)^3}\right)\frac{N^3}{3}-\frac{N^2}{2}\log N\log y+\frac{N^2}{4}(\log y)^2-\left(\gamma+\sum_p\frac{\log p}{p(p-1)}\right)\frac{N^2}{2}\log N\\&\quad-\left(\gamma-\frac{3}{2}+\log2\pi+\sum_p\frac{\log p}{p(p-1)^2}\right)\frac{N^2}{2}\log y+O(N^2)+O_\epsilon(N^{\frac32}y^{\frac12+\epsilon})+O(N^2\log(2N) y^{-\frac{1}{2}}).
\end{split}\end{equation}
\end{theorem}
\noindent{\it Proof.} The definition of $\mathfrak{S}_y$ and the formula $\displaystyle c_q(-k)=\sum_{\substack{d\vert q\\d\vert k}}d\!\;\mu\!\left(\frac{q}{d}\right)$ give
\begin{equation}\label{eq7.2} S_2=\sum_{q\leq y}\frac{\mu(q)^2}{\phi(q)^2}\sum_{d\vert q}d\!\;\mu\!\left(\frac{q}{d}\right)\sum_{\substack{1\leq k\leq N\\ d\vert k}}(N-k)^2\mathfrak{S}(k) .\end{equation}
Letting $k=dm$, the inner sum is
\begin{equation}\label{eq7.3} S_{21}(d)=d^2\sum_{1\leq m\leq N/d}(N/d-m)^2\mathfrak{S}(dm)=d^2S_{22}(N/d), \end{equation}
where
\begin{equation}\label{eq7.4} S_{22}(x)=\sum_{1\leq m\leq x}(x-m)^2\mathfrak{S}(dm)=2\int_1^x\sum_{1\leq m\leq u}(u-m)\mathfrak{S}(dm)\,du=:2\int_1^xS_{23}(u)\,du. \end{equation}

To evaluate this using Lemma 5, we write it in terms of $\mathfrak{G}_d$:
\begin{equation}\label{eq7.5}\begin{split}
S_{23}(x)&=\sum_{1\leq n\leq x}(x-n)\mathfrak{S}(dn)=\sum_{\substack{1\leq n\leq x\\2\vert dn}}(x-n)\frac{d}{(d,2)\phi(d)}\mathfrak{G}_d(n)\\&=\begin{cases}\displaystyle\frac{d}{2\phi(d)}\sum_{1\leq n\leq x}(x-n)\mathfrak{G}_d(n)& \mbox{if $d$ is even,}\\\displaystyle\frac{2d}{\phi(d)}\sum_{1\leq n\leq \frac{x}{2}}\left(\frac x2-n\right)\mathfrak{G}_d(n) & \mbox{if $d$ is  odd.}\end{cases}
\end{split}\end{equation}
The contribution to $S_{23}(x)$ from the main term of Lemma 5 is $\frac{d}{2\phi(d)}x^2$ regardless of the parity of $d$, and because $\displaystyle\sum_{p\vert 2d}\frac{\log p}{p-1}=\sum_{p\vert d}\frac{\log p}{p-1}$ if $d$ is even while $\displaystyle\log(x/2)+\sum_{p\vert2d}\frac{\log p}{p-1}=\log x+\sum_{p\vert d}\frac{\log p}{p-1}$ if $d$ is odd, the second term contributes $\displaystyle-\frac{x}{2}\Bigg(\!\!\!\;\log x+\gamma-1+\log2\pi+\sum_{p\vert d}\frac{\log p}{p-1}\:\!\Bigg)\!\!\;,$ again regardless of $d$'s parity.  The error term in Lemma 5 is $\ll_\epsilon x^{\frac{1}{2}}d^\epsilon$ and $\frac{d}{\phi(d)}\ll_\epsilon d^\epsilon$.  Thus
\begin{equation}\label{eq7.6} S_{23}(x)=\frac{d}{\phi(d)}\frac{x^2}{2}-\frac{x}{2}\Bigg(\!\!\!\;\log x+\gamma-1+\log 2\pi+\sum_{p\vert d}\frac{\log p}{p-1}\;\!\Bigg)\!+O_\epsilon(x^{\frac{1}{2}}d^\epsilon).\end{equation}
Integrating, and denoting $\gamma-\frac{3}{2}+\log 2\pi$ by $c_1$,
\begin{equation}\label{eq7.8}S_{22}(x)=2\int_1^xS_{23}(u)\,du=\frac{d}{\phi(d)}\frac{x^3}{3}-\frac{x^2}{2}\Bigg(\!\!\!\;\log x+c_1+\sum_{p\vert d}\frac{\log p}{p-1}\;\!\Bigg)-\frac{d}{3\phi(d)}+\frac{1}{2}\sum_{p\vert d}\frac{\log p}{p-1}+O_\epsilon(x^{\frac{3}{2}}d^\epsilon).\end{equation}
Thus, because $\frac{d}{\phi(d)}$ and $\sum_{p\vert d}\frac{\log p}{p-1}$ are both $O_\epsilon(d^\epsilon)$,
\begin{equation}\label{eq7.10}S_{21}(d)=\frac{N^3}{3\phi(d)}-\frac{N^2}{2}\Bigg(\!\!\!\;\log N-\log d+c_1+\sum_{p\vert d}\frac{\log p}{p-1}\;\!\Bigg)\!+O_\epsilon(N^{\frac{3}{2}}d^{\frac{1}{2}+\epsilon}).\end{equation}

For square-free $q$, $\displaystyle\sum_{d\vert q}\frac{d\mu\!\left(\frac{q}{d}\right)}{\phi(d)}=\mu(q)\prod_{p\vert q}\left(1-\frac{p}{p-1}\right)=\prod_{p\vert q}\frac{1}{p-1}=\frac{1}{\phi(q)}$, so the term $\dfrac{N^3}{3\phi(d)}$ contributes
\begin{equation}\label{eq7.11} \sum_{q\leq y}\frac{\mu(q)^2}{\phi(q)^2}\sum_{d\vert q}\frac{d\mu\!\left(\frac{q}{d}\right)\!N^3}{3\phi(d)}=\frac{N^3}{3}\sum_{q\leq y}\frac{\mu(q)^2}{\phi(q)^3}\end{equation}
to $S_2$.  Next, the terms $-\frac{N^2}{2}(\log N+c_1)$ are easily dealt with, and contribute
\begin{equation}\label{eq7.12}\begin{split}
-\frac{N^2}{2}(\log N+c_1)&\sum_{q\leq y}\frac{\mu(q)^2}{\phi(q)^2}\sum_{d\vert q}d\mu\!\left(\frac{q}{d}\right)=-\frac{N^2}{2}(\log N+c_1)\sum_{q\leq y}\frac{\mu(q)^2}{\phi(q)}\\&=-\frac{N^2}{2}(\log N+c_1)\left(\log y+\gamma+\sum_p\frac{\log p}{p(p-1)}\right)+O(N^2\log(2N) y^{-\frac{1}{2}})
\end{split}\end{equation}
by Lemma 3.  The error $O_\epsilon(N^{\frac{3}{2}}d^{\frac{1}{2}+\epsilon})$ contributes
\begin{equation}\label{eq7.12.5}\sum_{q\leq y}\frac{\mu(q)^2}{\phi(q)^2}\sum_{d\vert q}O_\epsilon(N^{\frac{3}{2}}d^{\frac{3}{2}+\epsilon})=O_\epsilon\!\Bigg(\!\!\!\;N^{\frac{3}{2}}\sum_{q\leq y}\frac{\mu(q)^2}{\phi(q)^2}q^{\frac{3}{2}+\epsilon}\;\!\Bigg)\!=O_\epsilon(N^{\frac32}y^{\frac12+\epsilon}).\end{equation}

For the remaining terms, we first evaluate the inner sum:
\begin{equation}\label{eq7.13}\begin{split}
\sum_{d\vert q}d\mu\!\left(\frac{q}{d}\right)\frac{N^2}{2}\Bigg(\!\!\!\;\log d-\sum_{p\vert d}\frac{\log p}{p-1}\;\!\Bigg)&=\frac{N^2}{2}\sum_{d\vert q}d\mu\!\left(\frac{q}{d}\right)\sum_{p\vert d}\left(1-\frac{1}{p-1}\right)\log p\\&=\frac{N^2}{2}\sum_{p\vert q}\Bigg(\frac{p-2}{p-1}\log p\sum_{\substack{d\vert q\\p\vert d}}d\mu\!\left(\frac{q}{d}\right)\!\!\!\;\Bigg)\!=\frac{N^2}{2}\sum_{p\vert q}\frac{p-2}{p-1}p\,\phi\!\left(\frac{q}{p}\right)\!\!\;\log p\\&=\frac{N^2}{2}\phi(q)\sum_{p\vert q}\frac{p(p-2)}{(p-1)^2}\log p=\frac{N^2}{2}\phi(q)\Bigg(\!\!\!\;\log q-\sum_{p\vert q}\frac{\log p}{(p-1)^2}\:\!\Bigg)
\end{split}\end{equation}
Thus the contribution of the terms $\frac{N^2}{2}\left(\log d-\sum_{p\vert d}\frac{\log p}{p-1}\right)$ is $\frac{N^2}{2}$ times
\begin{equation}\label{eq7.17}\sum_{q\leq y}\frac{\mu(q)^2}{\phi(q)}\Bigg(\!\!\!\;\log q-\sum_{p\vert q}\frac{\log p}{(p-1)^2}\:\!\Bigg)\!=\log y\sum_{q\leq y}\frac{\mu(q)^2}{\phi(q)}-\sum_{q\leq y}\frac{\mu(q)^2}{\phi(q)}\log(y/q)-\sum_{p\leq y}\frac{\log p}{(p-1)^2}\sum_{\substack{q\leq y\\p\vert q}}\frac{\mu(q)^2}{\phi(q)}\end{equation}
The first sum is evaluated in Lemma 3, and contributes
\begin{equation}\label{eq7.17.5}(\log y)^2+\left(\gamma+\sum_p\frac{\log p}{p(p-1)}\right)\log y+O(1).\end{equation}
Writing $q=pr$, the last sum is
\begin{equation}\label{eq7.18}\sum_{p\leq y}\frac{\log p}{(p-1)^3}\sum_{\substack{r\leq y/p\\(r,p)=1}}\frac{\mu(r)^2}{\phi(r)}=\sum_{p\leq y}\frac{\log p}{p(p-1)^2}(\log(y/p)+O(1))=\left(\sum_p\frac{\log p}{p(p-1)^2}\right)\log y+O(1).\end{equation}
We do the middle sum via contour integration:
\begin{equation}\label{eq7.19}\sum_{q\leq y}\frac{\mu(q)^2}{\phi(q)}\log(y/q)=\frac{1}{2\pi i}\int_{a-i\infty}^{a+i\infty} G(s)\frac{y^s}{s^2}ds,\quad\;\text{where}\quad\; G(s)=\sum_{n=1}^\infty\frac{\mu(n)^2}{\phi(n)}n^{-s}\;\;\text{and }a>0.\end{equation}
We have
\begin{equation}\label{eq7.20}G(s)=\prod_p\left(1+\frac{1}{(p-1)p^s}\right)=\zeta(s+1)\prod_p\left(1+\frac{1}{p(p-1)p^s}-\frac{1}{p(p-1)p^{2s}}\right)=:\zeta(s+1)H(s),\end{equation}
where $H(s)$ is analytic for $\re s>-1/2$ and $H(0)=1$.
Near $s=0$,
\begin{equation}\label{eq7.21}\begin{split}&G(s)\frac{y^s}{s^2}=\frac{1}{s^2}\zeta(s+1)H(s)y^s\\&=\frac{1}{s^2}\left(\frac{1}{s}+\gamma-\gamma_1s+O(s^2)\right)\left(1+H^\prime(0)s+\frac{H^{\prime\prime}(0)s^2}{2}+O(s^3)\right)\left(1+(\log y)s+\frac{(\log y)^2s^2}{2}+O(s^3)\right).\end{split}\end{equation}
The residue at 0 is then
\begin{equation}\label{eq7.22}\frac{(\log y)^2}{2}+(\gamma+H^\prime(0))\log y-\gamma_1+\gamma H^\prime(0)+\frac{H^{\prime\prime}(0)}{2},\end{equation}
so by moving the contour to $\re s=b$, $-1/2<b<0$ (convergence follows from \eqref{eq6.75}), we get
\begin{equation}\label{eq7.23}\sum_{q\leq y}\frac{\mu(q)^2}{\phi(q)}\log(y/q)=\frac{(\log y)^2}{2}+(\gamma+H^\prime(0))\log y+O(1).\end{equation}
We have $\displaystyle\frac{H^\prime(s)}{H(s)}=\sum_p\frac{(2p^{-2s}-p^{-s})\log p}{p(p-1)+p^{-s}-p^{-2s}}$, so $\displaystyle H^\prime(0)=\sum_p\frac{\log p}{p(p-1)}$.
Then combining \eqref{eq7.23} with \eqref{eq7.17.5} and \eqref{eq7.18},
\begin{equation}\label{eq7.25}\sum_{q\leq y}\frac{\mu(q)^2}{\phi(q)}\Bigg(\!\!\!\;\log q-\sum_{p\vert q}\frac{\log p}{(p-1)^2}\:\!\Bigg)=\frac{1}{2}(\log y)^2-\left(\sum_p\frac{\log p}{p(p-1)^2}\right)\log y+O(1).\end{equation}

Combining with the other terms \eqref{eq7.11}, \eqref{eq7.12} and \eqref{eq7.12.5},
\begin{equation}\label{eq7.26}\begin{split}
S_2&=\sum_{k\leq N}(N-k)^2\mathfrak{S}(k)\mathfrak{S}_y(k)\\&=\Bigg(\!\!\;\sum_{q\leq y}\frac{\mu(q)^2}{\phi(q)^3}\:\!\Bigg)\frac{N^3}{3}-\frac{N^2}{2}(\log N+c_1)\left(\log y+\gamma+\sum_p\frac{\log p}{p(p-1)}\right)+O(N^2\log(2N) y^{-\frac{1}{2}})\\&\quad+\frac{N^2}{2}\left(\frac{1}{2}(\log y)^2-\left(\sum_p\frac{\log p}{p(p-1)^2}\right)\log y+O(1)\right)+O_\epsilon(N^{\frac32}y^{\frac12+\epsilon})\\&=\Bigg(\!\!\;\sum_{q\leq y}\frac{\mu(q)^2}{\phi(q)^3}\:\!\Bigg)\frac{N^3}{3}-\frac{N^2}{2}\log N\log y+\frac{N^2}{4}(\log y)^2-\left(\gamma+\sum_p\frac{\log p}{p(p-1)}\right)\frac{N^2}{2}\log N\\&\quad-\left(c_1+\sum_p\frac{\log p}{p(p-1)^2}\right)\frac{N^2}{2}\log y+O(N^2)+O_\epsilon(N^{\frac32}y^{\frac12+\epsilon})+O(N^2\log(2N) y^{-\frac{1}{2}}),
\end{split}\end{equation}
as claimed.

\section{The sum $S_3$} In this section we prove the following result on $S_3$. 
 
\begin{theorem} We have, for $1\le  y\le \sqrt{N}$, 
\begin{equation}\begin{split}\label{eq8.1}S_3 :=\sum_{k\le N}\left(N-k\right)^2\mathfrak{S}_y(k)^2=\Bigg(\!\!\;\sum_{q\leq y}\frac{\mu(q)^2}{\phi(q)^3}\:\!\Bigg)&\frac{N^3}{3}-\left(\log y+\gamma+\sum_p\frac{\log p}{p(p-1)}\right)^2\frac{N^2}{2}\\ & +O(N^2y^{-\frac12}\log y)+O(Ny^2).\end{split}\end{equation}
\end{theorem}

\noindent{\it Proof.}
The definition of $\mathfrak{S}_y(k)$ and the formula $c_q(-k)=\displaystyle\sum_{\substack{d\mid q\\d\mid k}}d\,\mu\!\left(\frac{q}{d}\right)$ give
\begin{equation}\label{eq8.2}
S_3=\sum_{q\leq y}\sum_{q^\prime\leq y}\frac{\mu(q)^2}{\phi(q)^2}\frac{\mu(q^\prime)^2}{\phi(q^\prime)^2}\underbrace{\sum_{1\leq k\leq N}(N-k)^2c_q(-k)c_{q^\prime}(-k)}_{S_{31}},\end{equation}
and
\begin{equation}\label{eq8.3}S_{31}=\sum_{d\mid q}\sum_{d^\prime\mid q^\prime}d\,\mu\!\left(\frac{q}{d}\right)d^\prime\,\mu\Big(\frac{q^\prime}{d^\prime}\!\!\;\Big)\sum_{\substack{1\leq k\leq N\\ [d,d^\prime]\mid k}}(N-k)^2.\end{equation}
Using Lemma 1 on 
\[ \sum_{\substack{1\leq k\leq N\\ [d,d^\prime]\mid k}}(N-k)^2=[d,d^\prime]^2\sum_{\substack{1\leq k\leq \frac{N}{[d,d^\prime]}}}\left(\frac{N}{[d,d^\prime]}-k\right)^2,\]
we obtain
\begin{align*}
S_{31}&=\sum_{d\mid q}\sum_{d^\prime\mid q^\prime}[d,d^\prime]^2dd^\prime\mu\!\left(\frac{q}{d}\right)\mu\Big(\frac{q^\prime}{d^\prime}\!\!\;\Big)\left(\frac{N^3}{3[d,d^\prime]^3}-\frac{N^2}{2[d,d^\prime]^2}+O\!\left(\frac{N}{[d,d^\prime]}\right)\right)\\&=\frac{N^3}{3}\sum_{d\mid q}\sum_{d^\prime\mid q^\prime}(d,d^\prime)\mu\!\left(\frac{q}{d}\right)\mu\Big(\frac{q^\prime}{d^\prime}\!\!\;\Big)-\frac{N^2}{2}\sum_{d\mid q}\sum_{d^\prime\mid q^\prime}dd^\prime\mu\!\left(\frac{q}{d}\right)\mu\Big(\frac{q^\prime}{d^\prime}\!\!\;\Big)+O\!\!\;\Bigg(\!N \sum_{d\mid q}\sum_{d^\prime\mid q^\prime}d^2{d^\prime}^2 \:\!\Bigg)\!, 
\end{align*}
where we use that $\frac{N}{6[d,d^\prime]}+O(1)=O\!\left(\frac{N}{[d,d^\prime]}\right)$ because $[d,d^\prime]\leq dd^\prime\leq qq^\prime\leq y^2\leq N$.  
Thus
\begin{equation} \label{eq8.4} S_{3} = A_1(y)\frac{N^3}{3}-A_2(y)\frac{N^2}{2}+O(A_3(y)N) ,\end{equation}
where
\begin{equation}\label{eq8.5} A_1(y)=\sum_{q\leq y}\sum_{q^\prime\leq y}\frac{\mu(q)^2}{\phi(q)^2}\frac{\mu(q^\prime)^2}{\phi(q^\prime)^2}\sum_{d\mid q}\sum_{d^\prime\mid q^\prime}(d,d^\prime)\mu\!\left(\frac{q}{d}\right)\!\mu\!\left(\frac{q^\prime}{d^\prime}\right)\!,\end{equation}
\begin{equation}\label{eq8.6}A_2(y)=\sum_{q\leq y}\sum_{q^\prime\leq y}\frac{\mu(q)^2}{\phi(q)^2}\frac{\mu(q^\prime)^2}{\phi(q^\prime)^2}\Bigg(\!\!\;\sum_{d\mid q}d\mu\!\left(\frac{q}{d}\right)\!\!\!\;\Bigg)\!\!\:\Bigg(\!\!\;\sum_{d^\prime\mid q^\prime}d^\prime\mu\Big(\frac{q^\prime}{d^\prime}\!\!\;\Big)\!\!\!\;\Bigg)=\Bigg(\!\!\;\sum_{q\leq y}\frac{\mu(q)^2}{\phi(q)^2}\cdot\phi(q)\Bigg)^{\!\!\;\!2},\end{equation}
and
\begin{equation}\label{eq8.7}\begin{split}
A_3(y)&=\sum_{q\leq y}\sum_{q^\prime\leq y}\frac{\mu(q)^2}{\phi(q)^2}\frac{\mu(q^\prime)^2}{\phi(q^\prime)^2}\Bigg(\!\!\;\sum_{d\mid q}d^2\Bigg)\Bigg(\!\!\;\sum_{d^\prime\mid q^\prime}{d^\prime}^2\Bigg)=\Bigg(\!\!\;\sum_{q\leq y}\frac{\mu(q)^2}{\phi(q)^2}\sum_{d\mid q}d^2\!\;\Bigg)^{\!\!2}\\&=\Bigg(\!\!\:\sum_{dr\leq y}\frac{\mu(dr)^2}{\phi(dr)^2}d^2\!\:\Bigg)^{\!\!2}=\Bigg(\!\!\;\sum_{r\leq y}\frac{\mu(r)^2}{\phi(r)^2}\sum_{\substack{d\leq y/r\\(d,r)=1}}\frac{\mu(d)^2d^2}{\phi(d)^2}\!\;\Bigg)^{\!\!2}\\&\ll \Bigg(\!\!\;\sum_{r\leq y}\frac{\mu(r)^2}{\phi(r)^2}\cdot\frac{y}{r}\!\;\Bigg)^{\!\!2}\ll y^2,
\end{split}\end{equation}
using Lemma 2 for the last two steps.
We compute $A_1(y)$ the same way we did $B_1(y)$ in \S3, using \[ \displaystyle (d,d^\prime)=\sum_{\substack{r\mid d\\r\mid d^\prime}}\phi(r)\] to get
\begin{equation}\label{eq8.8} A_1(y)=\sum_{r\leq y}\phi(r)\left(\rule{0in}{.35in}\right.\!  \sum_{\substack{q\leq y\\r\mid q}}\frac{\mu(q)^2}{\phi(q)^2}\sum_{\substack{d\mid q\\r\mid d}}\mu\!\left(\frac{q}{d}\right)\!\!\!\!\;\left.\rule{0in}{.35in}\right)^2=\sum_{r\le y}\frac{\mu(r)^2}{\phi(r)^3}.\end{equation}
We conclude
\begin{equation}\label{eq8.10}S_3=\Bigg(\!\!\;\sum_{r\leq y}\frac{\mu(r)^2}{\phi(r)^3}\!\:\Bigg)\frac{N^3}{3}-\Bigg(\!\!\;\sum_{q\leq y}\frac{\mu(q)^2}{\phi(q)}\!\:\Bigg)^{\!\!2}\frac{N^2}{2}+O(Ny^2).\end{equation}
Theorem 6 now follows from Lemma 3.

\section{Completion of the proof of Theorem 2} By \eqref{eq5.1} and Theorems 4, 5, and 6, for $1\leq y\leq \sqrt{N}$ we have
\[\begin{split}\sum_{k\leq N}&(N-k)^2\widetilde{\mathfrak{S}}_y(k)^2=S_1-2S_2+S_3\\
&=\left(\prod_p\left(1+\frac{1}{(p-1)^3}\right)-2\sum_{q\leq y}\frac{\mu(q)^2}{\phi(q)^3}+\sum_{q\leq y}\frac{\mu(q)^2}{\phi(q)^3}\right)\frac{N^3}{3}\\
&\quad-\frac14N^2(\log N)^2+2\left(\frac12N^2\log N\log y-\frac14N^2(\log y)^2\right)-\frac12N^2(\log y)^2\\
&\quad+\left(\frac34 -\gamma -\frac12\log 2\pi -\frac12\sum_p\frac{\log p}{(p-1)^2}\right)N^2\log N\\
&\qquad+2\left(\left(\gamma+\sum_p\frac{\log p}{p(p-1)}\right)\frac{N^2}{2}\log N+\left(\gamma-\frac{3}{2}+\log2\pi+\sum_p\frac{\log p}{p(p-1)^2}\right)\frac{N^2}{2}\log y\right)\\
&\qquad-2\left(\gamma+\sum_p\frac{\log p}{p(p-1)}\right)\frac{N^2}{2}\log y\\
&\quad+O(N^2)+O(N^2\log(2N)y^{-\frac12})+O_\epsilon(N^{3/2}y^{1/2+\epsilon})+O(N^2y^{-\frac12}\log y)+O(Ny^2)\\
&=\left(\sum_{q>y}\frac{\mu(q)^2}{\phi(q)^3}\right)\frac{N^3}{3}+N^2\left(-\frac14(\log N)^2+\log N\log y-(\log y)^2\right)\\
&\quad+\left(\frac34-\frac12\log 2\pi+\sum_p\frac{(p-2)\log p}{2p(p-1)^2}\right)N^2\log N+\left(-\frac32+\log 2\pi+\sum_p\frac{(2-p)\log p}{p(p-1)^2}\right)N^2\log y\\
&\quad+O(N^2)+O(N^2\log(2N)y^{-\frac12})\\
&=\left(\sum_{q>y}\frac{\mu(q)^2}{\phi(q)^3}\right)\frac{N^3}{3}-\frac14N^2\left(\log\frac{N}{y^2}\right)^2+cN^2\log\frac{N}{y^2}+O(N^2)+O(N^2\log(2N)y^{-\frac12}).
\end{split}\]


\begin{thebibliography}{99}

\bibitem{Buttkewitz}  Y. Buttkewitz, {\it Exponential sums over primes and the prime twin problem,} Acta Math. Hungar., {\bf 131} (1-2) (2011), 46-58. DOI: 10.1007/s10474-010-0015-9 .

\bibitem{FriedlanderGoldston1995}J. B. Friedlander and D. A. Goldston, {\it Some singular series averages and the
distribution of Goldbach numbers in short intervals,} Illinois J. of Math.,
{\bf 39} No. 1, Spring 1995, 158--180.

\bibitem{Goldston1984} D. A. Goldston, {\it The second moment for prime numbers,} Quart. J. Math.
Oxford (2) {\bf 35} (1984), 153--163.

\bibitem{Goldston1990} D. A. Goldston {\it Linnik's theorem on Goldbach numbers in short intervals,}
Glasgow Math. J. {\bf 32} (1990), 285--297.

\bibitem{Goldston1992} D. A. Goldston, {\it On Hardy and Littlewood's contribution to the Goldbach conjecture,} in: E. Bombieri et al. (Eds.), {\it Proceedings of the Amalfi Conference on Analytic Number Theory} (1992), 115-155. 

\bibitem{Goldston2000} D. A. Goldston, {\it The major arcs approximation for  an exponential sum over primes,} Acta Arithmetica {\bf XCII.2} (2000), 169--179.

\bibitem{HardyLittlewood1922} G. H. Hardy and J. E. Littlewood, {\it Some problems of `Partitio numerorum'; III: On the expression of a number as a sum of primes}, Acta Math. {\bf 44} (1922), no. 1, 1--70. Reprinted as pp. 561--630 in {\it Collected Papers of G. H. Hardy}, Vol. I, Clarendon Press, Oxford University Press, Oxford, 1966.

\bibitem{Hildebrand1984} A. Hildebrand, \textit{ \"{U}ber die punktweise Konvergenz von Ramanujan-Entwicklungen
zahlentheoretischer Funktionen}, Acta Arithmetica \textbf{XLIV} (1984),  109--140.


\bibitem{Ingham1932}A. E. Ingham,  {\it The Distribution of Prime Numbers} (Cambridge Tracts in Mathematics and Mathematical Physics 30, Cambridge Univ. Press, Cambridge, 1932).

\bibitem{MontgomerySound1999} H.L. Montgomery and K. Soundararajan, {\it Beyond pair correlation},  Paul Erd\H{o}s and his mathematics, I (Budapest, 1999), 507--514, Bolyai Soc. Math. Stud., 11, Janos Bolyai Math. Soc., Budapest, 2002.
\
\bibitem{MontgomeryVaughan1973}H. L. Montgomery and R. C. Vaughan, {\it Error terms in additive prime number theory,} Quart. J. Math. Oxford (1), {\bf 24} (1973), 207--216.

\bibitem{MontgomeryVaughan2007} H. L. Montgomery and R. C. Vaughan, {\it Multiplicative Number Theory}, Cambridge Studies in Advanced Mathematics {\bf 97}, Cambridge University Press, Cambridge, 2007.

\bibitem{Selberg1947}  Atle Selberg, {\it On an elementary method in the theory of primes}. Norske Vid. Selsk. Forh. Trondheim {\bf 19} (1947) 64�67.

\bibitem{Titchmarsh} E. C. Titchmarsh, {\sl The Theory of the Riemann Zeta-Function,} 
2nd ed., revised by D. R. Heath-Brown, Clarendon (Oxford), 1986.

\bibitem{Vaughan2001} R. C. Vaughan, {\it On a variance associated with the distribution of primes in arithmetic progressions,} Proc. London Math. Soc. (3) {\bf 82} (2001), 533--553.

\bibitem{Wilson1922} B. M. Wilson, {\it Proofs of some formulae enunciated by Ramanujan}, Proc. London Math. Soc. (2), {\bf 21} (1922), 235--255. 

\bibitem{Zhang2014} Yitang Zhang, {\it Bounded gaps between primes,} Annals of Mathematics,  {\bf 179} (2014), (3),  1121--1174.

\end{thebibliography}
\end{document}